# Multiply Connected Topological Economics, Confidence Relation and Political Economy

*Yi-Fang Chang*
*Department of Physics, Yunnan University, Kunming 650091, China*
(e-mail: yifangchang1030@hotmail.com)

**Abstract**
Using the similar formulas of the preference relation and the utility function, we propose the confidence relations and the corresponding influence functions that represent various interacting strengths of different families, cliques and systems of organization. Since they can affect products, profit and prices, etc., in an economic system, and are usually independent of economic results, therefore, the system can produce a multiply connected topological economics. If the political economy is an economy chaperoned polity, it will produce consequentially a binary economy. When the changes of the product and the influence are independent one another, they may be a node or saddle point. When the influence function large enough achieves a certain threshold value, it will form a wormhole with loss of capital. Various powers produce usually the economic wormhole and various corruptions.

## 1. Introduction

In mathematical economics the fixed-point theorems of topology are used to prove the Nash equilibrium for n-person games. Arrow and Debreu presented a general model of Walrasian equilibrium theory, and proved the existence theorem of equilibrium for a competitive economy by topology [1]. Then McKenzie [2], Debreu [3], et al., developed the competitive equilibrium theory. Differential topology is introduced into economics, and Debreu [4] discussed two detailed questions.

Usual economic theories include only some interpretation on pure market process. The public choice theory is a great economics that intersects the two disciplines: the institutions are those of political science, and the method is that of economic theory [5-7]. It applies and develops scientific economic methods to other social regions. The public choice theory emphasizes comparative institutional analysis and, in particular, by their concentration on the necessary relationship between economic and political institutions. But alternative institutions may also have defects.

In the classical economics various quantities may be classified to two types: 1. Quantitative quantities, for example, capital, labor, product and profit, etc; 2. Qualitative quantities, for example, management, policy and preference, etc. The later possesses some human subjective factors. The two respects intersect usually each other, such the economic systems often show more complex social phenomena.

In the microeconomic theory of consumer behavior either a utility function or a binary relation can describe the preferences of an individual. The strict equivalence of these two primitive concepts, ordinal utility functions and preference relations, was first axiomatized by Debreu [8,9]. He studies the concept of cardinal utility in three different situations by means of the same mathematical result that gives a topological characterization of three families of parallel straight lines in a plane [8], and discussed that for every continuous complete and transitive binary relation $\geq$ defined on an arbitrary subset X of the commodity space R, there is a continuous utility representation; that is, there is a continuous function u of X into R such that $u(y) \geq u(x)$ if and only if $y \geq x$. Therefore, the more basic concept of preferences is applied instead of utility by means of a topology or a metric on the space of preferences. Undoubtedly, it is a great contribution for economics.

The topological structure on the space of preferences is very useful. For example, Hildenbrand used the structure to describe an exchange economy by its distribution of agents characteristics: preferences and endowments.

Moreover, the solution to the static maximization problem in Bellman equation yields a policy function that gives the optimal value of the current control as a function $g_t(x_t)$ of time and the current state. So the tomorrow state in is given by $x_{t+1} = m_t[g_t(x_t), x_t]$, and a solution to a similar problem then yields tomorrow's optimal control [10].

## 2. The confidence relations and the influence function

In microeconomics we introduce the confidence relations that represent various interacting strengths of different families, cliques and systems of organization. It is an important human relation in economics, even is independent of economic results. The confidence relation can be defined by a similar method with the preference relation in consumer theory.

The confidence relation $\geq$ defined on the choice set X is a complete preordering, continuous and strictly monotone. This requires [10]

(1). Reflexivity: $\forall x \in X, x \geq x$;

(2). Completeness: $\forall x, y \in X$, either $x \geq y$ or $y \geq x$ or both;

(3). Transitivity: $\forall x, y, z \in X$, $[x \geq y$ and $y \geq x] \Rightarrow x \geq z$.

Then $\geq$ can be represented by a real-valued, continuous and increasing payoff function.

Further, the definition of the influence function I is similar with the utility function: A real-valued function $I^i : X^i \to R$ represents a confidence preordering $\{\geq_i\}$ defined on the choice set $X^i$ of agent i if $\forall x, y \in X^i, x \geq_i y \Leftrightarrow I^i(x) \geq I^i(y)$. The influence function that represents a confidence preorder is not uniquely defined. Any monotonically increasing transformation $\varphi(\ )$ of I( ) will represent exactly the same confidences, because with $\varphi(\ )$ strictly increasing, we have

$$I(x) \geq I(y) \text{ if and only if } \varphi[I(x)] \geq \varphi[I(y)], \qquad (1)$$

for all $\forall x, y \in X$. Hence I( ) is an ordinal influence function. The sign of the difference I(x)-I(y) is important because it tells us which outcome is confided, but the value of this difference is meaningless, as it will change with any nontrivial increasing transformation $\varphi(\ )$. It is also a basic characteristic of topology, where those concrete spacing values are meaningless. Although the influence function is similar to the utility function that obeys the law of diminishing marginal utility, but the influence function seems to obey the law of augmenting lust for power.

## 3. Multiply connected topological economy

The confidence relation, the corresponding influence function I( ) and the function $\varphi(\ )$ can affect products Q, profit and prices, etc., in an economic system. But, they are usually independent of economic results, and sometimes are stochastic, even change suddenly. In a continuous topological manifold of economics they break easily original structure, and form a new hole or branch region. This will construct a multiply connected topological manifold. In an image the economic structure is a cup, while the influence function is a handle.

In a multiply connected region of topology there is a famous Euler-Poincare formula

$$\sum_{m=1}^{n}(-1)^m a_m = \sum_{m=1}^{n}(-1)^m p_m . \qquad (2)$$

For a convex polyhedron, $a_0, a_1, a_2$ denote the number of vertices, edges, and faces, respectively; $p_m$ is the mth Betti number of complex K. This may be considered intuitively as the numbers of m-dimensional holes in K, or is the number of (m+1)-dimensional chains that must be added to K so that every free m-cycle on K is a boundary [11]. The number $\sum_{m=1}^{n}(-1)^m a_m$ is called the Euler characteristic of the complex K. In the polyhedron $p_0 = p_2 = 1, p_1 = 2p$, p is the deficiency of a curved surface. In 2-dimensional curved surface, $a_0 = a_1 + 1 - 2p$. Assume that vertices represent the number of market, which is direct proportional to the sales volume *y* and the profit, and edges represent the market network. But the multiply connected economy brings the profit decrease. In this case there is a defective profit due to the deficiency *p*.

In some systems of organization the profit maximization and the confidence relations are inseparable. The aim of a pure producer is the profit maximization

$$\pi(y,w) = \max_{x,y}\{py - wx\}, \quad (3)$$

where *y* and *x* are output and input, *p* and *w* are output and input prices. For a social system with the influence function I, we should define an aim function as

$$A = \pi + I. \quad (4)$$

The objective function $f(x;\alpha)$ gives his payoff, when he faces environment $\alpha$ and chooses action *x* [10]. Usual profit is

$$\pi(Q) = TR(Q) - TC(Q), \quad (5)$$

where *TR(Q)* and *TC(Q)* are the total revenue and the total cost, the three quantities all are the functions of products *Q=f(K,L)*. Now the maximum principle is the aim function maximization. The initial condition of the profit maximization is

$$\frac{d\pi}{dQ} = \frac{dTR}{dQ} - \frac{dTC}{dQ} = 0, \quad (6)$$

$$\therefore \frac{dTR}{dQ} = \frac{dTC}{dQ}, \quad (7)$$

i.e., *MR=MC*, the marginal revenue equals to the marginal cost. While now we derive a new result:

$$\frac{dA}{dQ} = \frac{dTR}{dQ} - \frac{dTC}{dQ} + \frac{dI}{dQ} = 0, \quad (8)$$

$$\therefore MR - MC = -\frac{dI}{dQ} < 0. \quad (9)$$

In the social system the aim function deviates the profit maximization, which is defective usually since $MR - MC < 0$.

The production function in the traditional theory of the firm expresses output Q as a function of two inputs: capital K and labor L,

$$Q=Q(K,L). \quad (10)$$

The profit maximization of the firm is

$$\pi = pQ - mK - nL, \quad (11)$$

where p, m and n are the prices of output, capital and labor flows respectively [12]. If the influence function regards as a condition of the economic system, the economic meaning of the influence function will also be able to be discussed using the Lagrange method of conditional maximization.

**4. Political economy and multiply connected topological economy**

A complete market economy should be a simplex free economy. But, a non-market economy, and any oligarch economy must be a multiply connected topological economy, which changes to a higher dimension, for example, it will add a new dimension with man-rule. This economy may be dismembered and comminuted, and has various holes and be mangled easily for distortions.

According to <Oxford Advanced Learner Dictionary of Current English> the political economy is study of the political problems of government. It as an early title is very famous, for example, David Ricardo's <The Principles of Political Economy and Taxation>(1817), Thomas Robert Malthus' <Principles of Political Economy>(1820) and <Definitions of Political Economy>(1827), James Mill's <Elements of Political Economy>(1821), and Karl Marx's <The Critique of Political Economy> and so on. The political economy now sounds old-fashioned but usefully emphasizes the importance of choice between alternatives in economics which remains, despite continuing scientific progress [13].

If the political economy is an economy chaperoned polity, it will produce consequentially a binary economy. Its basic group is different with a complete market economy in the algebraic topology.

A general change of the supply-demand function is

$$\frac{dQ_d(Q_s)}{dt} = f(Q_d, Q_s, p) + V + S. \tag{12}$$

Here V is a governmental potential, and S is a stochastic factor. The equation (12) has the outside force and the potential V, such the economic results change along with different V. If V is inequitable and factitious, the results will possess bigger stochasticity.

Various powers are some different attractors, and produce the economic wormhole and various corruptions. In particular, if the multifarious confidence relations exist, a whole economic system and corresponding topological manifold will be covered with many big and small holes like bruises and scars. In this society the highest economic aim is only the confidence relation for families, cliques and systems of organization. Some big or small powers cling to an economic system, and form the multiply connected topological economy, and have a series of corruption with the self-similarity. It is a special type of the fractal economy. There is a binary economic function of power-business. Both is usually asymmetry, i.e., is inequality. Under this system various aspects tend spontaneously to the breaking of symmetry. An imperium is an economic black hole, which will derive the huge corruption, and finally the system dies out.

The political economy should be a pair coupling equations on polity and economy. Assume that a potential is $U = 2aX^2Y$. Here X is a confidence relation, Y is an economic benefit, and *a* is a coefficient. From this the difference between a theoretic value and a practice value will be estimated.

When the politics is put in command, the economy and its equation will be neglected. But, when oligarch notices the social crisis, the economic rules will be obeyed more. Both alternation exhibits a periodicity. If they conflict, the result will be reform for the economics bigger than politics, or be retrogression for the economics smaller than politics.

The political economy is usually imperfect economic question, even completely is not an economic question for some particular cases. It is not a strict economic rule, because in this case economy is only an appendage of polity. The economy will change along with polity.

The multiply connected topological economy may be extended to various relations between economy and other politics, family, religion, etc. Further, it may be developed to many regions of without direct relations with economy, for example, welfare, environment, and full employment, etc.

**5. Singularity and wormhole**

If the influence function changes as time, the system will be more complex. Assume that the economic system and its change are linear [10]:

$$\frac{dQ}{dt} = a_{11}Q + a_{12}I, \tag{13}$$

$$\frac{dI}{dt} = a_{21}Q + a_{22}I. \tag{14}$$

Their characteristic matrix is

$$\begin{pmatrix} a_{11} & a_{12} \\ a_{21} & a_{22} \end{pmatrix}. \tag{15}$$

The corresponding characteristic equation is

$$\lambda^2 - (a_{11} + a_{22})\lambda + (a_{11}a_{22} - a_{12}a_{21}) = \lambda^2 - T\lambda + D = 0. \tag{16}$$

From this we may discuss the general cases. As the simplest example, if the changes of the product and the influence are independent one another, i.e., $a_{12} = a_{21} = 0$, the solutions of the equations will be $Q = Q_0 \exp(a_{11}t), I = I_0 \exp(a_{22}t)$. $\Delta = T^2 - 4D = (a_{11} - a_{22})^2 \geq 0$ for a real domain.

When $a_{11}, a_{22}$ are real numbers of the same signs, D>0, the state $(Q_0, I_0)$ of the system is a node point, which is stable for $a_{11}, a_{22} <0$, and is unstable for $a_{11}, a_{22} >0$ (Fig.1).

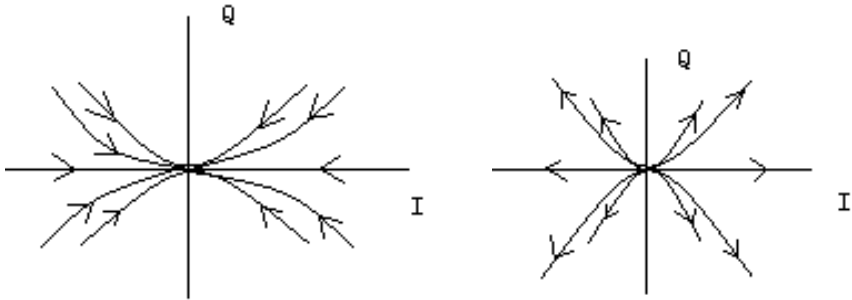

Fig.1   Stable and unstable node points

When $a_{11}, a_{22}$ are real numbers of opposite signs, $D<0$, the state $(Q_0, I_0)$ of the system is a saddle point (Fig.2).

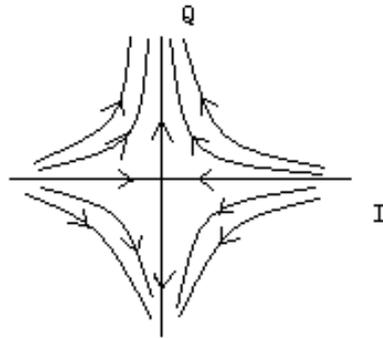

Fig.2   Saddle point

It represents that product increases and the conference decreases. If the two changes of the product and the influence intersect one another, the states of the economic system will be able also to be the spiral (focal) point, or center, etc.

The form of the influence function can be an unrestricted function, even a stochastic function. Perfect competition prevails that each producer and consumer regards the prices paid and received as independent of his own choices [1]. An economy with the confidence relations and the influence functions is a type of imperfect competitive economic systems, and break the symmetries in economic topology. They are not homeomorphic spaces. Usually this structure will hinder the economic development. If the confidence relations and the influence functions have $p$-levels or $p$-types, i.e., $\sum_{i=1}^{p} I_i(Q)$, they will construct a multiply connected normal curved surface with the deficiency $p$.

When the influence function large enough achieves a certain threshold value, the economic elasticity of topological structure will be broken, and a new hole will appear. Unified market economy will be riddled with holes. This will form a new multiply connected topological manifold. As an example, using the concept of general relativity a large influence as mass of general relativity forms a pit in the economic system. According to Fuller-Wheeler theory [14], a very strong pit can construct a wormhole, sometimes called the Einstein-Rosen bridge [15]. Therefore, some capital will pass through a throat into another topological space, or from a region to another region in the same space (Fig.3). This model will may describe a loss of capital (including waste, and corruption).

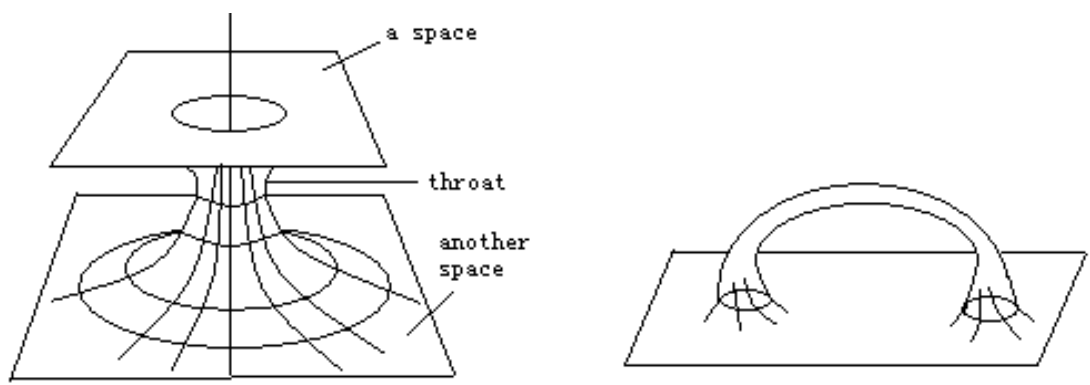

　　Fig.3 The wormhole model in social economics from a space into another topological space, or from a region to another region in the same space

　　In a word, the confidence relation and the influence function provide the useful tools for a description of human activity in economic system. Based on the main characteristics of a new epoch of knowledge economy and its similarity with the information theory, we proposed the four theorems of the knowledge economic theory, and expounded the production function and basic equations, and discussed some possible directions of the development on the knowledge economy and its theory [16]. We think that topology and its tools in economy will exhibit larger function in the epoch.

　　Finally, this method of the multiply connected topological economy can be extended to various aspects on polity-law, on polity-eduction, on government-people and so on.